\documentclass[11pt]{article}
\usepackage{amsmath, amssymb, epsfig, graphicx}
\usepackage{bm}
\usepackage{mathrsfs}
\usepackage{color}
\usepackage{fullpage} 
\usepackage[small,bf]{caption}
\newtheorem{theorem}{Theorem}[section]

\newtheorem{corollary}[theorem]{Corollary}
\newtheorem{lemma}[theorem]{Lemma}

\newtheorem{remark}{Remark}

\numberwithin{equation}{section}
\newlength{\algorithmwidth}
\newcommand{\bigO}{\mathrm{O}}

\newcommand{\notate}[1]{#1}
\newcommand{\defby}{\overset{\mathrm{\scriptscriptstyle{def}}}{=}}

\def \P {\mathbb{P}}
\def \E {\mathbb{E}}
\def \< {\langle}
\def \> {\rangle}

\DeclareMathOperator*{\argmax}{arg max} 
\def \endprf{\hfill {\vrule height6pt width6pt depth0pt}\medskip}
\newcommand{\pr}[2]{\langle {#1} , {#2} \rangle}

\usepackage{fullpage}

\begin{document}
\bibliographystyle{plain}

\pagestyle{plain}

\title{Acceleration of Randomized Kaczmarz Method\\via the Johnson-Lindenstrauss Lemma}
\author{Yonina C. Eldar$^1$, Deanna Needell$^2$\thanks{Corresponding author: Deanna Needell. Email: dneedell@stanford.edu}\\
  \vspace{-.1cm}\\
  $^1$Department of Statistics, Stanford University, Stanford, CA 94305\\
    \vspace{-.3cm}\\
  $^2$Department of Electrical Engineering, Technion - Israel Institute of Technology, Haifa 32000}
  %\thanks{Partially supported by NSF DMS EMSW21-VIGRE grant and the Israel Science Foundation under Grant no. 1081/07.}
  
   %\keywords{Kaczmarz method, randomized Kaczmarz method, computer tomography, signal processing.}
\date{\today}

\maketitle

\vspace{-0.3in}

\begin{abstract}
The Kaczmarz method is an algorithm for finding the solution to an overdetermined \notate{consistent} system of linear equations $Ax=b$ by iteratively projecting onto the solution spaces.  The randomized version put forth by Strohmer and Vershynin yields provably exponential convergence in expectation, which for highly overdetermined systems even outperforms the conjugate gradient method.  In this article we present a modified version of the randomized Kaczmarz method which at each iteration selects the optimal projection from a randomly chosen set, which in most cases significantly improves the convergence rate.  We utilize a Johnson-Lindenstrauss dimension reduction technique to keep the runtime on the same order as the original randomized version, adding only extra preprocessing time.  We present a series of empirical studies which demonstrate the remarkable acceleration in convergence to the solution using this modified approach.
\end{abstract}

%\subjclass{68W20, 65T50, 41A46}

\section{Introduction}

The Kaczmarz method~\cite{K37:Angena} is a popular algorithm for solving overdetermined \notate{consistent} systems of linear equations.  Due to its simplicity and speed, it has been used in a variety of applications ranging from tomography to digital signal processing~\cite{CFMSS92,SS87:Applications,N86:TheMath}.  The method uses a series of alternating projections to iteratively converge to the solution of $Ax=b$, and is therefore computationally feasible even for very large systems.  Given an initial guess $x_0$ and denoting by $a_1, \ldots, a_m$ the rows of the $m \times n$ matrix $A$, each iteration of the method orthogonally projects the current estimation onto the next hyperplane defined as the solutions to $\< a_i, x \> = b_i$, chosen in a cyclic fashion.  The algorithm can be described by the iterations:
$$
x_{k+1} = x_k + \frac{b[i] - \pr{a_i}{x_k}}{\|a_i\|_2^2}a_i,
$$
where $x_k$ is the $k^{th}$ iterate, $b[i]$ (here and throughout) denotes the $i$th coordinate of $b$, and $i = (k$ mod $m) + 1$. 

Although this technique has been used in practice for quite some time, theoretical guarantees on convergence were difficult to obtain~\cite{DH97:Therate,G05:Onthe,HN90:Onthe}.  It is clear that by design of the algorithm, the convergence rate depends on the ordering of the rows in $A$. Therefore, poorly ordered rows can lead to slower convergence.  To overcome this difficulty, the rows of $A$ can be selected in a random fashion.  It has been observed that this randomized version of the algorithm improves the convergence rate~\cite{N86:TheMath,HM93:Algebraic}, however only recently have theoretical results been obtained~\cite{SV06:Arandom,SV09:Arand,N10:rknoisy}.

\subsection{Randomized Kaczmarz} In~\cite{SV06:Arandom,SV09:Arand}, Strohmer and Vershynin propose at each iteration to randomly select a row of $A$ with probability proportional to the Euclidean norm of the row.  The randomized Kaczmarz (RK) method can thus be described by
\begin{equation}\label{eq:rk}
x_{k+1} = x_k + \frac{b[p(i)] - \pr{a_{p(i)}}{x_k}}{\|a_{p(i)}\|_2^2}a_{p(i)},
\end{equation}
where $p(i)$ takes values in $\{1, \ldots, m\}$ with probabilities $\frac{\|a_{p(i)}\|_2^2}{\|A\|_F^2}$.  Here and throughout, $\|A\|_F$ denotes the Frobenius norm of $A$ and $\|\cdot\|_2$ denotes the standard Euclidean norm or spectral norm for vectors or matrices, respectively.  \notate{The selection rule used here is not optimal in general. The motivation for setting the rule according to the weight of the row norms is two-fold.  First, it allows for a guarantee of expected exponential convergence for the Kaczmarz method~\cite{SV09:Arand}. Second, it is a computationally efficient strategy since often these values will be known approximately or exactly, and will only need to be computed once.  A selection rule of this type is of course also related to the idea of preconditioning the matrix $A$ by scaling its rows.  Although other diagonal preconditioners may certainly perform better in general, finding such an optimal preconditioner is itself an optimization problem of high complexity.  With the above selection strategy,} the following exponential bound was shown in~\cite{SV06:Arandom,SV09:Arand} for the convergence in expectation of this randomized method:
\begin{equation}\label{SVbound}
\mathbb{E}\|x_k - x\|_2^2 \leq \Big(1 - \frac{1}{R}\Big)^k\|x_0 - x\|_2^2,
\end{equation}
where $R = \|A^{-1}\|^2\|A\|_F^2$ and $x_0$ is an arbitrary initial estimate.  Since we will always assume that $A$ has full column rank, the norm $\|A^{-1}\| \defby \inf\{M : M\|Ax\|_2 \geq \|x\|_2$ for all $x\}$ is well-defined.  This bound is essentially independent of the number of rows of $A$.  Moreover, the bound shows that for well conditioned matrices $A$, the RK method yields expected exponential convergence to the solution in just $\bigO(n)$ iterations (see Section 2.1 of~\cite{SV09:Arand}).  Since each iteration consists of a single projection taking $\bigO(n)$ time, this shows that the overall method has $\bigO(n^2)$ runtime, which is clearly superior to other methods such as Gaussian elimination which takes $\bigO(mn^2)$, especially when the system is very large.  The discussion in~\cite{SV09:Arand} shows that the randomized Kaczmarz method often even outperforms the celebrated conjugate gradient method.  For example, when $A$ is a Gaussian matrix and $m > 3n$, the RK method provably requires fewer computations, and empirical studies show that this improvement is substantial.  See Section 4.2 of~\cite{SV09:Arand} for details. 

The empirical and theoretical benefits of this approach lead one to ask whether it is also accurate in the more realistic case when noise is present.  One may thus consider the (now possibly inconsistent) system $Ax \approx b + w$ where $w$ is an arbitrary error vector that has been added to the consistent system $Ax = b$.  It is shown in~\cite{N10:rknoisy} that in this case we have exponential convergence to the solution within an error factor:
$$
\mathbb{E}\|x_k - x\|_2 \leq \Big(1 - \frac{1}{R}\Big)^{k/2}\|x_{0}\|_2 + \sqrt{R}\gamma,
$$
where $R$ is the same as above and $\gamma = \max_i \frac{|w[i]|}{\|a_i\|_2}$.  It is also shown that this bound is sharp and is attained even for simple examples~\cite{N10:rknoisy}.   

\subsection{Modified approach} 
To further improve the convergence rate of the RK method, we suggest a different approach to selecting the rows of $A$.  \notate{Although our ideas should also apply seamlessly to the case when noise is present and the system becomes inconsistent, in this work we only consider the noiseless case and leave a detailed analysis in the presence of noise for future work.}  Since the projections in the algorithm~\eqref{eq:rk} are orthogonal, it can be seen that the optimal projection in the $k$th iteration is the one that maximizes $\|x_{k+1} - x_k\|_2$.  By definition of the iterations~\eqref{eq:rk}, one can calculate these quantities by computing inner products between the rows $a_i$ of $A$ and the current iterate $x_k$.  Since computing one inner product requires $\bigO(n)$ operations, we clearly cannot afford to perform more than a constant number of these in a given iteration.  Our approach, therefore, is to project the rows of $A$ onto a lower dimensional space in such a way that the geometry of the vectors is approximately preserved.  We then perform calculations of the form~\eqref{eq:rk} with respect to these low dimensional vectors, and select the best projection.

By construction, our modified algorithm will converge to the solution of $Ax=b$ \textit{in the worst case} as fast as the standard RK method.  In practice, we expect the convergence to be much faster, especially when the lower dimension $d$ onto which we project the rows is not too small.  The improvement in each iteration can be quantified in terms of $d$ and the current estimation $x_k$, as we demonstrate in Section~\ref{sec:just}.  In Section~\ref{sec:sims} we demonstrate that empirically our technique outperforms the standard RK method in terms of convergence rate.  The runtime of this modified algorithm of course depends on the dimension $d$ onto which we project the rows of $A$.  If $d \ll n$, then each iteration will require $\bigO(dn)$ operations, meaning that the overall runtime for expected exponential convergence becomes \textit{at most} $\bigO(dn^2)$.  There is a tradeoff in the choice of the dimension $d$.  If $d$ is small, then the runtime per iteration remains small.  However, if $d$ is too small then our technique reduces to the standard RK method, and we will not gain in the convergence rate.  In the next section we discuss the runtime and implementation, and show why the selection of $d$ on the order of $\log n$ is the right choice.  This gives a \textit{worst case} runtime of $\bigO(n^2\log n)$, although we expect a much faster convergence as we also see in simulations.  This worst case runtime however, is still the same as that of the standard RK method up to the log factor.

%The proven exponential convergence rates and the abundance of applications of the method lead one to ask whether the convergence can be accelerated even further.  
\section{Implementation and Runtime}\label{sec:imp}
Since the projections in the algorithm are orthogonal, one easily sees that the optimal projection in the $k$th iteration would be the one that maximizes $\|x_{k+1} - x_k\|_2$, or equivalently, the term
\begin{equation}\label{eq:term}
\frac{|b[i] - \pr{a_{i}}{x_k}|}{\|a_{i}\|_2}.
\end{equation}
Unfortunately, calculating this term takes $\bigO(n)$ time, so that to keep the overall runtime at $\bigO(n^2)$ one can only afford to make this computation a constant number of times.  However, if we could significantly reduce the dimension of the vectors $a_i$ and $x_k$ used in the calcuation, then more of these calculations could be done at each iteration, and the best out of those computed could be chosen, leading to accelerated convergence.  Our idea is thus to project the vectors onto a low dimensional space such that the geometry is preserved.  This will allow approximation of the inner products $\< a_i, x_k \> $ and the norms $\|a_i\|_2$ (if they are not known a priori) from the projected data.  To do this, we will consider a Johnson-Lindenstrauss type projection.  The well-known Johnson-Lindenstrauss Lemma~\cite{JL84:Extensions} states that with high probability, there is a projection of a finite set of points onto a space logarithmic in the number of points that approximately preserves geometry.  This can be summarized as follows.

\begin{lemma}[Johnson-Lindenstrauss~\cite{JL84:Extensions}]\label{lem:JL} Let $\delta > 0$ and let $S$ be a finite set of points in $\mathbb{R}^n$.  Then for any $d$ satisfying 
\begin{equation}\label{eq:k}
d \geq C\frac{\log |S|}{\delta^2},
\end{equation}
there exists a Lipschitz mapping $\Phi : \mathbb{R}^n \rightarrow \mathbb{R}^d $ such that
\begin{equation}\label{eq:JL}
(1-\delta)\|s_i - s_j\|_2^2 \leq \|\Phi(s_i) - \Phi(s_j)\|_2^2 \leq (1+\delta)\|s_i - s_j\|_2^2,
\end{equation}
for all $s_i, s_j\in S$, where $C$ is an absolute constant.
\end{lemma}

\notate{\textbf{Remark.}  The value of $C$ in which this lemma holds depends on the distribution from which $\Phi$ is created.  When $\Phi$ is Gaussian, one has $C\leq 8$~\cite{dasgupta2003elementary}.}

Although this lemma as stated only guarantees \textit{existence} of such a mapping, in their proof the map $\Phi$ is chosen as the projection onto a random $d$-dimensional subspace of $\mathbb{R}^n$.  This result has been improved over time and now one can easily construct such a (random) projection which preserves the geometry (see e.g.~\cite{Ach03:Database-Friendly-Random,DG99:elem,IM98:Approx}).  Indeed, it is shown in~\cite{Ach03:Database-Friendly-Random} that whenever a distribution satisfies certain moment conditions, the $d\times n$ random matrix $\Phi$ whose entries are chosen i.i.d. with respect to that distribution will satisfy~\eqref{eq:JL} with high probability provided $d$ satisfies~\eqref{eq:k}.  The Guassian distribution, for example, satisfies these moment conditions; therefore the matrix with i.i.d. Gaussian entries will preserve geometry with high probability.  Recently there has been work on constructing transforms which satisfy~\eqref{eq:JL} but that also provide a fast multiply (see e.g.~\cite{AC09:fast,HV09:johnson}).  For example, Ailon and Chazelle construct in~\cite{AC09:fast} a $C\delta^{-2}\log|S| \times n$ transform $\Phi$ satisfying~\eqref{eq:JL} with high probability whose multiply requires roughly $n\log n + \delta^{-2}\log^3 |S|$ operations.  Hinrichs and Vybiral provide a multiply using $n\log n$ operations when $\Phi$ has slightly more rows, on the order of $\log^2|S|$.   Even more recently, Ailon and Liberty show in~\cite{AL10:Almost} that when the $d\times n$ $\Phi$ is the composition of a randomly subsampled Hadamard (or Fourier) matrix and a random sign matrix, then $\Phi$ satisfies~\eqref{eq:JL} with high probability when $d = C\delta^{-4}\log|S|\log^4 n$.  This matrix has an $n\log n$ multiply, and so this result provides an optimal fast JL transform, up to the power $-4$ on $\delta$ and the polylogarithmic dependence on $n$. 

\subsection{Implementation}

In our setting, the Johnson-Lindenstrauss Lemma allows us to project the rows of $A$ as well as the estimations $x_k$ onto a space of substantially lower dimension.  This will then let us approximately calculate the terms in~\eqref{eq:term} using far fewer operations, which we can use to decide on which hyperplane to project the current estimate.  %That is, if approximating the term~\eqref{eq:term} will now take $\bigO(d)$ time, we can now afford to make $n/d$ of these calculations.  
We note that the projection of the rows of $A$ will be performed offline, adding to the preprocessing time, whereas the projection of the estimation $x_k$ will be done at each iteration.  \notate{We choose to use an analagous strategy as in the RK algorithm for our row selection; that is, using the weight of the row norms.  This choice yields a worst case convergence rate which is the same as that guaranteed by the original RK method (see Remark 2 below).}
This leads to the following modified randomized Kaczmarz method, called Randomized Kaczmarz via Johnson-Lindenstrauss (RKJL) which can be summarized as follows.

\bigskip

\textsc{Randomized Kaczmarz via Johnson-Lindenstrauss (RKJL)}

\nopagebreak

\fbox{\parbox{\algorithmwidth}{
  \textsc{Input:} $m\times n$ matrix $A$, coefficient vector $b\in\mathbb{R}^m$, parameter $d$, initial estimate $x_0$
  
  \textsc{Output:} Approximate $x$ solving $Ax=b$

%  \textsc{Procedure:}
  \begin{description}
    \item[Initialize:] Set $k=0$, create a $d\times n$ Gaussian matrix $\Phi$ and set $\alpha_i = \Phi a_i$.  Repeat the following $\bigO(n)$ times:
    \item[Select:] Select $n$ rows so that each row $a_i$ is chosen with probability $\|a_i\|_2^2 / \|A\|_F^2$ as in~\ref{eq:rk}.  For each row selected, calculate 
    $$\gamma_i = \frac{|b[i] - \pr{\alpha_{i}}{\Phi x_k}|}{\|\alpha_{i}\|_2},$$
    and set $j = \argmax_i \gamma_i$. 
    \item[Test:] For $a_j$ and the first row $a_l$ selected out of the $n$, explicitly calculate
    $$
    {\gamma_j^*} = \frac{|b[j] - \pr{a_{j}}{x_k}|}{\|a_{j}\|_2} \quad\text{and}\quad {\gamma_l^*} = \frac{|b[l] - \pr{a_{l}}{x_k}|}{\|a_{l}\|_2}.
    $$
    If ${\gamma_l^*} > {\gamma_j^*}$, set $j=l$.
    \item[Project:] Set $$x_{k+1} = x_k + \frac{b[j] - \pr{a_j}{x_k}}{\|a_j\|_2^2}a_j.$$
    \item[Update:] Set $k = k+1$.
  \end{description}
 }}

    \bigskip

\begin{remarks}
{\bfseries 1.} We show in Section~\ref{sec:just} that $d = \bigO(\log n)$ will be enough to approximately preserve geometry and thus give convergence improvements.  Of course greater values of $d$ may give greater improvements on convergence, at the expense of more computational cost at each iteration.
%Since the randomized Kaczmarz method converges exponentially fast (in expectation) to the solution in $\bigO(n)$ iterations, we only plan on calculating $\bigO(n^2)$ inner products in the algorithm.  Thus we only ask that the Johnson-Lindenstrauss mapping $\Phi$ preserve the geometry on a set of roughly $n^2$ vectors, which is why we choose $d = C\delta^{-2}\log(n)$.  One could instead set choose $d$ larger, which would increase the probability of improved convergence when using the RKJL multiple times for multiple problem instances with the matrix $A$. 

{\bfseries 2.} In the Test stage of the algorithm, we see that in addition to approximating the $n$ inner products, we also exactly calculate the inner product of a randomly selected row and also the row that was chosen.  This will guarantee that the convergence is not slowed by any drastic consequences of the error in the approximations, and does not affect the overall runtime.  

{\bfseries 3.} We note that the initialization step may of course be computationally expensive, as is calculating the probabilities $p(i)$ in the standard version~\eqref{eq:rk}.  However, this need only be done once, and thus this version of the algorithm will be beneficial for situations in which the same matrix $A$ is used in many problems.  This is the case for many applications such as the wave-scattering problem~\cite{SPM89:CG} and structural mechanics problems~\cite{FR93:Implicit}; see~\cite{CW97:Ana,S00:ALevin} for others.
\end{remarks}

We next turn to an analysis of this modified method. In Section~\ref{sec:sims} we provide numerical results demonstrating the improved convergence rate.

%\section{Discussion}

%The RKJL algorithm allows one to speed up convergence to solutions while keeping rougly the same online computational cost, and only a one-time extra cost initially.  For applications like\notate{list them here}, the same measurement matrix $A$ is used repeatedly, and so the one-time cost of the initialization step is worth the accelerated convergence in subsequent runs of the algorithm.  The actual online runtime of the algorithm is the same as the standard randomized Kaczmarz method, up to a log factor, as we discuss next.

\subsection{Runtime}
As discussed above and in~\cite{SV06:Arandom,SV09:Arand}, the standard randomized Kaczmarz method converges exponentially fast to the solution in $\bigO(n)$ iterations, resulting in a total runtime of $\bigO(n^2)$.  The RKJL algorithm will thus also converge (in expectation) in at most $\bigO(n)$ iterations, and so it remains to calculate the runtime of each iteration.  

The first step in an iteration is to calculate $\Phi x_k$.  Since $\Phi$ is a $d \times n$ matrix, this computation in general takes $\bigO(nd)$ time.  Next, each $\alpha_i$ lives in a $d$ dimensional space, so that calculating inner products in the selection step costs only $\bigO(d)$.  Since we calculate $n$ such inner products, the total calculation time is $\bigO(nd)$.  %Thus when a fast transform is used, the overall runtime for the selection step is $\bigO(n\log n)$.  
The projection and update steps clearly take $\bigO(n)$ and $\bigO(1)$ time, respectively, leading to an overall runtime per iteration of $\bigO(nd)$.  Therefore, after $\bigO(n)$ iterations, we see that the overall runtime of the algorithm is $\bigO(n^2d)$.  Lemma~\ref{lem:noise} below shows that $d$ can be chosen on the order of $\log n$.  Therefore, RKJL converges exponentially fast in at most $\bigO(n^2\log n)$ time, which is the same as the runtime of standard RK, up to the log factor.  In practice, the runtime is much faster as we show in Section~\ref{sec:sims}.

Because the algorithm will need to use roughly $n^2$ rows of $A$ (assuming $n^2 \leq m$), the matrix $\Phi$ will have to be applied to at least this many vectors, resulting in a $\bigO(n^3d)$ initial cost.  If the algorithm is used repeatedly for various problem instances over the same matrix $A$, then one may wish to apply $\Phi$ to \textit{all} the rows of $A$, yielding a $\bigO(mnd)$ cost.  Even when $d=\bigO(\log n)$ this is of course substantial, but for applications in which the algorithm will be used many times, this one-time cost will become minimal.  

This initial computational cost occurs in other methods as well, for example, in submatrix selection algorithms that take $\bigO(mn^2)$ to select a well represented $n\log n\times n$ submatrix of $A$ (see e.g.~\cite{DMS06:Sampling,DMS08:Relative}).  These methods randomly select the submatrix (according to a particular probability model), and if the submatrix selected represents $A$ well, then it can be used to solve the system $Ax=b$.  However, there is some probability that the subsystem cannot be solved, in which case it must be reselected again, and so on.  This is in contrast to RKJL, for which we are always guaranteed convergence, and most likely with improved expected convergence rate.  The choice of method will of course depend on the application, and in some cases it may even be beneficial to use some combination of these approaches.  

\section{Analytical Justification}\label{sec:just}

We next analyze how the Johnson-Lindenstrauss Lemma is utilized by our method.  We will assume here that the system is real-valued and homogeneous (ie. $Ax=0$), that the rows of $A$ all have unit norm, and that the initial guess $x_0$ also satisfies $\|x_0\|_2 \leq 1$.  These assumptions are of course not necessary, but will make the analysis simpler.  \notate{We discuss the case where the row norms may be far from equal in Remark 3 below.}  We begin with an easy lemma which shows that the geometry of the vectors used in the RKJL method is approximately preserved. 

\begin{lemma}\label{lem:noise}
Let $\Phi$ be the $n\times d$ (Gaussian) matrix with $d = C\delta^{-2}\log(n)$ as in the RKJL method.  Set $\gamma_i = \< \Phi a_i, \Phi x_k\> $ also as in the method.  Then $|\gamma_i -  \<  a_i,  x_k\> | \leq 2\delta $ for all $i$ and $k$ in the first $\bigO(n)$ iterations of RKJL. 
\end{lemma}
\begin{remark} This lemma shows that with $d$ chosen on the order of $\bigO(\log n)$, the geometry of the vectors involved in the RKJL method is approximately preserved.  In practice, $d$ should thus be chosen of this order to gain improvements in convergence.
\end{remark}
\textit{Proof.}
We employ the Johnson-Lindenstrauss Lemma (Lemma~\ref{lem:JL}) with $\mathcal{S}$ consisting of all $a_i$ and $x_k$ in the first $\bigO(n)$ iterations of RKJL.  Then since $|\mathcal{S}| \lesssim n^2$, the condition~\eqref{eq:k} is satisfied, and so~\eqref{eq:JL} holds for all $a_i$ and $x_k$ used in the algorithm. By this and the parallelogram law, we have
\begin{align*}
\gamma_i &= \< \Phi a_i, \Phi x_k\> \\
&= \frac{1}{4}\Big(\|\Phi a_i + \Phi x_k\|_2^2 - \|\Phi a_i - \Phi x_k\|_2^2\Big) \\
&\leq \frac{1}{4}\Big((1+\delta)\| a_i +  x_k\|_2^2 - (1-\delta)\| a_i -  x_k\|_2^2\Big)\\
&= \< a_i, x_k \> + \frac{1}{4}\delta( \| a_i +  x_k\|_2^2 + \| a_i -  x_k\|_2^2 )\\
&\leq \< a_i, x_k \> + 2\delta.
\end{align*}

Similarly we have that $\gamma_i \geq \< a_i, x_k \> - 2\delta $, which completes the claim.
\endprf

This shows that the terms $\gamma_i$ used for selection in the algorithm are approximately equal to the actual desired values $\< a_i, x_k \> $.  Thus for $\delta$ small, the RKJL algorithm makes well educated decisions at each iteration which allows for quicker convergence (see Theorem~\ref{thm:real} below).  This also shows that when the estimation $x_k$ becomes very close to the true solution $x=0$, the error $\delta$ begins to dominate and improvements may no longer be expected.  However, this does not pose a problem since it only occurs when the estimate is already approximately $x$.

It is clear from construction of the RKJL algorithm (and especially in light of Remark 2 above), that convergence using RKJL is at least as fast as the standard randomized version.  Moreover, when the error produced by applying $\Phi$ is small, the RKJL method will project onto the ``best'' hyperplane out of those it selected in that iteration.  \notate{Since the probability of choosing this ``best'' row when selecting only a single row is strictly less than the probability of choosing that row when a set of rows is selected,} this implies that the only case in which RKJL would not provide a strictly faster convergence rate is when $\< a_i, x_k \> = \< a_j, x_k \> $ for all rows $a_i$, $a_j$ selected in the $k$th iteration.  

Given a current estimate $x_k$, one can explicitly compare the expectation of the improvement the next estimation provides, for both the RKJL and standard randomized methods.  First observe that if $P$ denotes the projection in the $k$th iteration, then $x_{k-1} - x_k$ resides in the kernel of $P$, and is thus orthogonal to the space onto which $P$ projects.  This space contains $x_k - x$ since $x$ is the solution to all equations $Ax=b$ and $x_k = Px_{k-1}$.  Therefore, $x_k - x$ and $x_{k-1} - x_k$ are orthogonal, implying that
\begin{equation}\label{eq:geom}
\|x_k - x_{k+1}\|_2^2 = \|x - x_{k}\|_2^2 - \|x - x_{k+1}\|_2^2.
\end{equation}
The relation~\eqref{eq:geom} shows that the larger $\|x_k - x_{k-1}\|_2$, the bigger the improvement made in that iteration.  We thus fix an estimation $x_{k}$ and analyze the expectation of $\|x_k - x_{k+1}\|_2$ for the RKJL method versus the standard one.  %For simplicity, we first perform these calculations under the assumption that the projection $\Phi$ exactly preserves geometry, meaning that $\gamma_j$ in the algorithm is precisely $x_{k+1} - x_k$.  This will give us good intution on what we can hope for in improvements using RKJL.  
For convenience, we again consider the real and homegenous case (ie. when $b = 0$), and assume the rows of $A$ have unit norms.  We then have the following result.

\begin{theorem}\label{thm:real}
 Fix an estimation $x_k$ and denote by $x_{k+1}$ and $x_{k+1}^*$ the next estimations using the RKJL and the standard RK method, respectively.  Set $\gamma_j^* = |\pr{a_j}{x_k}|^2$ and reorder these so that $\gamma_1^* \geq \gamma_2^* \geq \ldots \geq \gamma_m^*$.  Then when $d=C\delta^{-2}\log n$,
$$
\E \|x_{k+1} - x\|_2^2  \leq \min \left[\E\|x_{k+1}^*-x\|_2^2 - \sum_{j=1}^m \Big(p_j - \frac{1}{m}\Big)\gamma_j^* + 2\delta,\quad \E\|x_{k+1}^*-x\|_2^2\right],
$$
where 
$$
p_j = \Biggr\{ \begin{array}{ll} \frac{\binom{m-j}{n-1}}{\binom{m}{n}}, & j \leq m-n+1 \\ 0, & j > m-n+1 \end{array}
$$ 
are non-negative values satisfying $\sum_{j=1}^m p_j = 1$ and $p_1 \geq p_2 \geq \ldots \geq p_m = 0$. 
\end{theorem}

\textit{Proof.}
%The proof is similar to the proof of Lemma~\ref{lem:exact}, except 
Since we assume the rows of $A$ have unit norm and that $b=0$, we see that $\gamma_J$ is precisely the value of $\|x_{k+1} - x_k\|_2^2$ if the algorithm were to select row $J$.  We begin by examining the $k$th RKJL iteration.  

Let $\mathcal{L}$ denote the set of rows chosen in the selection step, so that $|\mathcal{L}|= n$. If exact geometry were preserved (i.e. $\delta=0$), then the method would simply select the index of the largest $\lambda_i^*$ contained in $\mathcal{L}$.  However, due to the error induced by $\Phi$, even when the ``best'' row is selected to be in $\mathcal{L}$, the algorithm may not choose this row for the projection.  

We thus define sets $T_j$ for $j=1,\ldots m$, which consist of rows which could be ``confused'' in this way,
$$
T_j \defby \{i: |\gamma_i^* - \gamma_j^*| \leq 2\delta\}, \quad\text{and}\quad \mu_j = \min\{\gamma_i^*: i\in T_j\}. 
$$
From Lemma~\ref{lem:noise}, if $j\in\mathcal{L}$ and $1,\ldots, j-1\notin\mathcal{L}$, then the worst row we could choose is one that would give $\|x_{k+1} - x_k\|_2^2 = \mu_j$.  Therefore,
\begin{align*}
\E \|x_{k+1} - x_k\|_2^2 &= \E_J \gamma_J^* \\
&\geq \mu_1\P(1\in \mathcal{L}) + \mu_2\P(2\in \mathcal{L}\cap 1\notin \mathcal{L}) + \ldots + \mu_{m-n+1} \P(1,2,\ldots,m-n\notin \mathcal{L})\\
&= \sum_{j=1}^{m-n+1} \mu_j \P(j\in \mathcal{L}\cap 1, \ldots, j-1\notin \mathcal{L}).
\end{align*}

Since each row of $A$ has equal norm, each $row$ of $A$ is equally likely to be selected in $\mathcal{L}$.  Thus 
$$
\P(j\in \mathcal{L}\cap 1, \ldots, j-1\notin \mathcal{L}) = \frac{\binom{m-j}{n-1}}{\binom{m}{n}}.
$$
Therefore, we have
$$
\E \|x_{k+1} - x_k\|_2^2 \geq \sum_{j=1}^{m-n+1} \mu_j p_j.
$$
%where we have set $p_j = \frac{\binom{m-j}{n-1}}{\binom{m}{n}}$ (and $p_j = 0$ for $j > m-n+1$).  Note that clearly we have $p_1 \geq p_2 \geq \ldots \geq p_m = 0$.  

Lemma~\ref{lem:noise} and the definition of $T_j$ guarantee that $\mu_j \geq \max(\gamma_j^* - 2\delta, 0)$.  This along with the fact that $\sum p_j = 1$ yields
\begin{equation}\label{eq:final}
\E \|x_{k+1} - x_k\|_2^2 \geq \Big(\sum_{j=1}^{m-n+1} \gamma_j^* p_j\Big) - 2\delta.
\end{equation}
Finally, since all the rows of $A$ have the same norm, the standard RK method selects each row uniformly at random, so that
\begin{equation}\label{eq:al}
\E \|x_{k+1}^* - x_k\|_2^2 = \sum_{j=1}^m \frac{1}{m}\gamma_j^*.
\end{equation}      
Combining~\eqref{eq:al} with~\eqref{eq:final} and~\eqref{eq:geom} we have
\begin{align*}
\E \|x_{k+1} - x\|_2^2 &= \E\|x_k - x\|_2^2 - \|x_k - x_{k+1}\|_2^2\\
&\leq \E\|x_k - x\|_2^2 - \Big(\sum_{j=1}^{m-n+1} \gamma_j^* p_j\Big) + 2\delta\\
&= \E \|x_{k+1}^* - x\|_2^2 + \E \|x_{k+1}^* - x_k\|_2^2  - \Big(\sum_{j=1}^{m-n+1} \gamma_j^* p_j\Big) + 2\delta\\
&= \E \|x_{k+1}^* - x\|_2^2 - \sum_{j=1}^m \Big(p_j - \frac{1}{m}\Big)\gamma_j^* + 2\delta.
\end{align*}
This along with the fact that by construction of RKJL, $\E \|x_{k+1} - x\|_2^2 \leq \E \|x_{k+1}^* - x\|_2^2$, completes the claim.
\endprf

\begin{remarks} 
{\bfseries 1.} Theorem~\ref{thm:real} gives a lower bound, which shows improvements in the ``worst case'', when the error induced by the Johnson-Lindenstrauss projection causes the method to choose a row of $A$ in the worst way.  Numerical experiments (as seen in the next section) demonstrate substantial improvements in the convergence rate.

{\bfseries 2.}
Note that since the sequences $\{\gamma_j^*\}$ and $\{p_j\}$ are non-increasing and $\sum_{j=1}^m p_j = 1 = \sum_{j=1}^m \frac{1}{m}$, the sum $\beta = \sum_{j=1}^m \Big(p_j - \frac{1}{m}\Big)\gamma_j^*$ is non-negative.    %$\E \|x_{k+1} - x_k\|_2^2 \geq \E \|x_{k+1}^* - x_k\|_2^2$.  
Furthermore, $\beta=0$ only when $\gamma_1^* = \gamma_2^* = \ldots = \gamma_m^*$.  Indeed, if $\gamma_i^* > \gamma_j^*$ even for just one pair $i < j$, then for $\delta$ small enough, the RKJL method provides strict convergence improvement.  Knowledge of $x_k$ and $A$ would allow one to precisely calculate the improvement.

{\bfseries 3. }The theorem is proven under the assumption that the rows of $A$ have the same norm.  The same argument holds (with different values of $p_j$) %.  To generalize, one would let $q_j$ be the probability that row $a_j$ corresponding to $\gamma_j^*$ is selected when only one row is selected (i.e. $q_j = \|a_j\|^2/\|A\|_F^2$), and $q_j^* = \P(j\in \mathcal{L}\cap 1, \ldots, j-1\notin \mathcal{L})$.  Then following the same argument as above one would have%  
still showing improvement when this assumption does not hold, and numerical experiments show similar results in either case.  \notate{Although the analysis of exact improvement in these other cases may quickly become quite complicated, we recall that by construction the RKJL method offers overall improvement in any case.} 

\notate{It is also helpful to identify some particularly interesting scenarios, for example, when   %No matter how the row norms are distributed, one follows the argument of Theorem~\ref{} to see that
%$$
%\E \|x_{k+1} - x\|_2^2  \leq \min \left[\E\|x_{k+1}^*-x\|_2^2 - \sum_{j=1}^m (q_j^*-q_j) \gamma_j^* + 2\delta,\quad \E\|x_{k+1}^*-x\|_2^2\right].
%$$
%where $q_j$.  This shows improvement in convergence under any scenario, however, the distribution of the row norms will of course affect the magnitude of such an improvement.  
 one or a few rows are substantially of larger norm.  In this case the standard RK algorithm (noticing that each row is selected independently of the previous selections) will choose this row repeatedly with high probability.  Clearly this may slow convergence (especially when these rows are highly correlated), and although there is still guaranteed exponential convergence, $R$ in this case is much larger so that the guaranteed rate is slower.  In RKJL, although the guaranteed worst case rate is the same, these large rows are likely not to be selected when their contributions toward the solution is minimal.  This will again speed up convergence, since in the language of Theorem~\ref{thm:real}, the $\lambda^*_k$ corresponding to the rows which are highly correlated with the previous projection will be such that $k \approx m$.  This means that the same argument as in Theorem~\ref{thm:real} will hold for all $\lambda^*_j$ with $j\leq k\approx m$.  Although this means the improvement in RKJL may not be as substantial as in the case of equal normed rows, we will still see a large improvement. } 
 
\notate{These ideas highlight the fact that the selection strategy is not optimal in general.  For example, if there are $k \ll m$ rows which are highly uncorrelated but have very small norm, and $m-k\approx m$ equal rows with substantially larger norm, neither RK nor RKJL will perform well.  Of course if this were reversed, with the uncorrelated rows having large norm and the equal rows having small norm, the selection strategy will yield excellent convergence in both RK and especially RKJL.  We emphasize again that the choice of the selection rule is certainly not optimal in general, but is computationally efficient and allows for provable expected exponential rate of convergence.  Choosing an optimal selection strategy for a given system is itself a problem of high complexity.  If one is using RKJL and investing preprocessing time to perform dimension reduction, one may also wish to simply normalize the rows to avoid such difficulties.} 

\end{remarks}

Theorem~\ref{thm:real} implies the following corollary, showing the improved convergence using RKJL when exact geometry is preserved (i.e. when $\delta\rightarrow 0$).

\begin{corollary}\label{cor:exact}
Fix an estimation $x_k$ and denote by $x_{k+1}$ and $x_{k+1}^*$ the next estimations using the RKJL and the standard method, respectively.  Set $\gamma_j^* = |\pr{a_j}{x_k}|^2$ and reorder these so that $\gamma_1^* \geq \gamma_2^* \geq \ldots \geq \gamma_m^*$.  Then when exact geometry is preserved ($\delta \rightarrow 0$),
$$
\E\|x_{k+1} - x\|_2^2 \leq \E\|x_{k+1}^*-x\|_2^2 - \sum_{j=1}^m \Big(p_j - \frac{1}{m}\Big)\gamma_j^*.
$$
\end{corollary}

%We can combine Corollary~\ref{cor:exact} and Lemma~\ref{lem:noise} to get a bound on the improvement using RKJL in the realistic setting when exact geometry is not preserved.  

\section{Numerical Results}\label{sec:sims}

\begin{figure}[ht]
\begin{center}
% $\begin{array}{c@{\hspace{.1in}}c}
\includegraphics[width=3in]{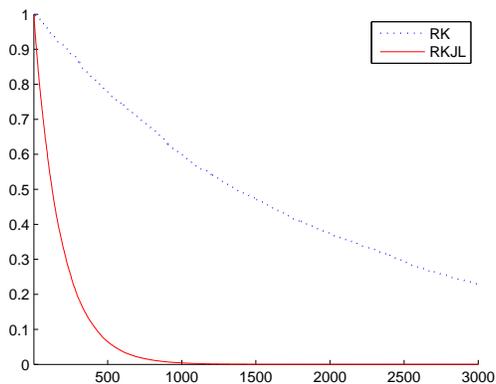}  %&
%\includegraphics[width=3in]{rk2f.eps}\\
%(a) & (b)
%\end{array}$
\end{center}
\caption{$\ell_2$-Error (y-axis) as a function of the iterations (x-axis). The dashed line is standard Randomized Kaczmarz, and the solid line is the modified one, without a Johnson-Lindenstrauss projection.  Instead, the best move out of the randomly chosen $n$ rows is used.  Note that we cannot afford to do this computationally.}\label{fig:fig1}
\end{figure}

We now demonstrate improved convergence using the RKJL method.  The first experiment we run is in the computationally infeasible situation where we do not use the Johnson-Lindenstrauss projection, but simply choose the best row out of the randomly selected $n$ rows.  This experiment will demonstrate the improved convergence using RKJL when the error $\delta$ induced by $\Phi$ goes to $0$. This is the best improvement one can hope for in RKJL.  When the problem sizes grow very large, the effect of the $\delta^{2}$ term in~\eqref{eq:k} becomes minimal, so it may be realistic to take $\delta$ quite small. For these simulations we use a $60000\times 1000$ matrix with Bernoulli entries and use a homogeneous system with an initial estimate chosen uniformly at random on the sphere.  We see in Figure~\ref{fig:fig1} that the convergence in this scenario is significantly improved. 

%The second plot in Figure~\ref{fig:fig1} depicts the results of a similar simulation, except that the Johnson-Lindenstrauss is used with $d = 800$ (and $\Phi$ is Gaussian).  This would be a realistic scenario, except that we use a different matrix $\Phi$ at each iteration.  We again cannot afford to do this computationally, however, it is interesting that the results are quite similar to those where no projection was used.

\begin{figure}[ht]
\begin{center}

\includegraphics[width=3in]{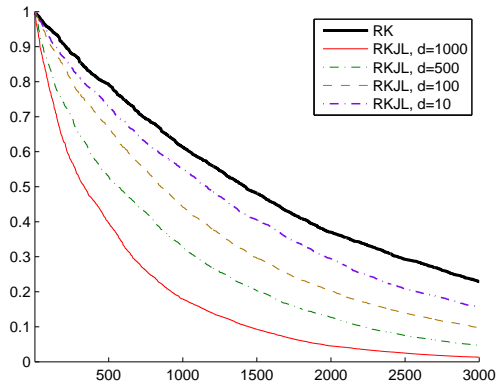}  
\end{center}
\caption{$\ell_2$-Error (y-axis) as a function of the iterations (x-axis) for various values of $d$ with $m=60000$ and $n=1000$.}\label{fig:fig2}
\end{figure}

Next we run simulations using a Johnson-Lindenstrauss matrix $\Phi$.  We generate $\Phi$ and the system $Ax=b$ the same as above, but now run RKJL using various values of $d$.  In Figure~\ref{fig:fig2} we see exactly what we expect, that with higher values of $d$ (corresponding to lower $\delta$ values), we have much quicker convergence.  The speedup in convergence using larger $d$ needs to be weighed against the increase in computation per iteration, as was discussed in Section~\ref{sec:imp}.  Since right now there are no theoretical guarantees on precisely how the convergence is affected by larger $d$, this comparison should be done empirically.  Finally, it is clear that as $m$ and $n$ grow large, the impact on $d$ of forcing $\delta$ to be small becomes minimal.  Thus for very large systems, using $d = \bigO(\log n)$ will give convergence that looks more like that in Figure~\ref{fig:fig1}.

%    \begin{figure}[ht]
%\begin{center}
% $\begin{array}{c@{\hspace{.1in}}c}
%\includegraphics[width=3in]{m_60000_1000.eps}  &
%\includegraphics[width=3in]{m_60000_800.eps}\\
%(a) & (b)\\
%\includegraphics[width=3in]{m_60000.eps}  &
%\includegraphics[width=3in]{m_60000_300.eps}\\
%(c) & (d)
%\end{array}$
%\end{center}
%\caption{$\ell_2$-Error (y-axis) as a function of the iterations (x-axis). The dashed line is standard RK, and solid is RKJL.  In all plots, $m=60000$ and $n=1000$. The values for $d$ are: (a) 1000, (b) 800, (c) 500, and (d) 300.}\label{fig:fig2}
%\end{figure} 
%
%    \begin{figure}[ht]
%\begin{center}
%
%\includegraphics[width=3in]{m_60000_all.eps}  
%\end{center}
%\caption{$\ell_2$-Error (y-axis) as a function of the iterations (x-axis) for various values of $d$ with $m=60000$ and $n=1000$.}\label{fig:fig3}
%\end{figure} 

%\notate{Not sure which figure looks best out of Figures 2, 3, and 4.}

%Our final experiment is to test the length of computation time for RK and RKJL to improve the estimation by a specified factor.  For these experiments we use the same setup but with $m=5000$, $n=1000$, and various levels of $d$.  We run both algorithms until they have improved the estimation by a factor of $2$, that is, $\|x_k - x\|_2 \leq 0.5\|x_0 - x\|_2$, and plot the length of computation time to do so.  While the code used is not necessarily completely optimized, both algorithms were implemented the same way, so that comparisons will be fair.  Figure~\ref{fig:times} shows the computation time for both algorithms.  We see that even with the additional runtime needed to multiply $\Phi$ and the estimation at each iteration, RKJL still outperforms RK significantly when $d=500$ or $d=1000$.
%
%    \begin{figure}[ht]
%\begin{center}
%
%\includegraphics[width=3in]{alltimes.eps}  
%\end{center}
%\caption{Computational time in seconds (y-axis) over several trials.  In all experiments, $m=60000$ and $n=1000$.  The mean time for RK was 2414 seconds, for RKJL with $d=100$ was 36.7 seconds, with $d=500$ was 53.3 seconds, and RKJL with $d=1000$ was $69.0$ seconds.}\label{fig:times}
%\end{figure} 

\subsection*{Acknowledgements}

This work is partially supported by the NSF DMS EMSW21-VIGRE grant and the Israel Science Foundation under Grant no. 1081/07.  We would also like to thank Emmanuel Cand\`es and Thomas Strohmer for helpful suggestions.

\small
\bibliography{rk}

\end{document}